\documentclass{article}
\usepackage{graphicx} 
\usepackage[T2A, T1]{fontenc}
\usepackage[utf8]{inputenc}
\usepackage{amsmath}
\usepackage{amssymb}
\usepackage{amsthm}
\usepackage{amsfonts}
\usepackage{mathtext}
\usepackage{geometry}
\usepackage{dsfont}
\usepackage[unicode, pdftex]{hyperref}
\title{Criticality conditions in the Derrida–Retaux model with a random
number of terms}
\author{Alexey Lotnikov, Anna Kotova }
\date{2023}
\newtheorem{lemma}{Lemma}
\newtheorem{theorem}{Theorem}

\begin{document}
 \maketitle
\section*{Abstract}
The article considers the Derrida--Retaux model with a random number of terms, i.e. a sequence of integer random variables
defined by the relations
$ X_{n + 1} = (X_n^{(1)} + X_n^{(2)} + ... + X_n^{(N_n)} - a)^{+}$, $n\ge 0$,
where $X_n^{j}$ are independent copies of $X_n$, the values of $N_j$ are independent and identically distributed, $a$ is a positive integer.
The energy in the model is defined as $Q:=\lim\limits_{n\to\infty} \frac{\mathbb{E}(X_{n})}{(\mathbb{E}N_1)^{n}}$. We present sufficient
conditions (in terms of distributions of $X_0$ and $N_1$) for subcritical ($Q=0$) and supercritical ($Q>0$) regimes of model behavior.
\section{Introduction}
\subsection{Derrida --- Retaux model}

The theory of fluctuations in chaotic systems is a branch of physics that studies the behavior of systems in which chaos and randomness are present. The purpose of the theory is to understand how fluctuations or small changes in the properties of a system can lead to global changes in its overall behavior.

In particular, two states may be distinguished in the system and transitions between them are studied. The Derrida --- Retaux model\cite{derrit} in a simplified form describes such transitions (depinning transition) between certain states called "pinned" and "unpinned", which we do not need to describe. Similar relationships have also been studied in the context of physical \cite{wetting} and mathematical studies \cite{collet}.

The Derrida --- Retaux model with a random number of terms, which is studied in this paper, is formulated as follows. The model parameters are an integer $a > 0$ and two integer random variables: $X_0\ge 0$ – the initial value, and $N\ge 1$ – the number of terms. It is assumed that $X_0$ has a finite first moment and is not a constant, as well as $\mathbb{P}(N > 1) > 0$. The functioning of the system is determined by the recurrence relation:
\begin{eqnarray}
\label{1}
 X_{n + 1} = (X_n^{(1)} + X_n^{(2)} + ... + X_n^{(N_{n+1})} - a)^{+},
\end{eqnarray}
where $X_n^{(1)}$ and $X_n^{(2)}$  are independent copies of $X_n$, and $ \forall \ x \in \mathbb{R}$, $x^{+} := \max(x, 0)$  is
the positive part of $x$.

\noindent Let us note that $$\mathbb{E}(\sum\limits_{j=1}^{N_{n+1}} X_n^{(j)} - a) \leq \mathbb{E}(\sum\limits_{j=1}^{N_{n+1}} X_n^{(j)} - a)^{+} \leq \mathbb{E}\sum\limits_{j=1}^{N_{n+1}} X_n^{(j)}.$$

\noindent Also $\mathbb{E}\sum\limits_{j=1}^{N_{n+1}} X_n^{(j)} = \mathbb{E} N_{n+1} \cdot \mathbb{E}X_n = \mathbb{E} N \cdot \mathbb{E}X_n$, therefore
$$\mathbb{E}N \cdot \mathbb{E} X_n - a \leq \mathbb{E} X_{n+1} \leq \mathbb{E}N \cdot \mathbb{E} X_n,$$
$$\cfrac{\mathbb{E}(X_{n+1})}{(\mathbb{E}N)^{n+1}} \leq \cfrac{\mathbb{E}(X_{n})}{(\mathbb{E}N)^{n}} \quad\text{    and    }\quad \cfrac{\mathbb{E}(X_n) - \cfrac{a}{\mathbb{E}N-1}}{(\mathbb{E}N)^n} \leq \cfrac{\mathbb{E}(X_{n+1}) - \cfrac{a}{\mathbb{E}N-1}}{(\mathbb{E}N)^{n+1}}.$$
Thus, the following limit is well-defined
\begin{eqnarray}
\label{2}
    {Q} = \lim\limits_{n \rightarrow \infty} \downarrow \cfrac{\mathbb{E}(X_{n})}{(\mathbb{E}N)^{n}} = \lim\limits_{n \rightarrow \infty} \uparrow \cfrac{\mathbb{E}(X_n) - \cfrac{a}{\mathbb{E}N-1}}{(\mathbb{E}N)^n}.
\end{eqnarray}

The parameter $Q$ is commonly referred to as free energy, and the recurrence relation described above was first introduced in the article by Derrida and Retaux \cite{derrit} and serves to describe the Depinning transition process.
Studying this process is important in both mathematics and physics. Of particular interest is the dependence of the energy value $Q$ on the initial data --- the random variable $X_0$ and $a$. The parameter $a$ is referred to as the tax.

Depending on the value of $Q$, two cases are distinguished: supercritical ($Q$ > 0) and subcritical ($Q$ = 0). A critical case is also distinguished, which essentially <<separetes>> the situations described above. The most interesting questions regarding this model arise in critical cases or cases close to them.

The recurrence relation \hyperref[1]{(1)} with parameter $a$ = 1 and determined numbers of terms was considered in the article \cite{lifsh}. This article provides necessary and sufficient conditions on the random variable $X_0$ under which supercritical, subcritical, and critical cases occur, respectively.

We attempted to generalize the results obtained in this article to arbitrary $a \geq 1$ and random numbers of terms. We found sufficient conditions on the random variable $X_0$ guaranteeing subcriticality or supercriticality.

\section{Main result}
The main result of our work is the following theorem.\\ 
Let $F_0, G$ be the moment generating function of $X_0$ and $N$ respectively.

\begin{theorem}
Let it be $D_0(s, m) = (m - 1)sF'_0(s) - aF_0(s)$.\\
1) If $D_0(\mathbb{E}N^{\frac{1}{a}}, \mathbb{E}N) > 0, \text{ then } Q > 0$.\\
2) Let $\exists M : \mathbb{P}( N \leq M ) = 1$. If $D_0(1 + \frac{M - 1}{a}, M) <0, \text{ then } Q = 0.$
\end{theorem}

\noindent Let us emphasize that the first point is true without assuming that the random variable $N$ is bounded.

\noindent The result of Theorem 1 is interesting even if the number of terms is fixed, i.e. $\mathbb{P}(N = n) = 1$ for some $n$. In this case, Theorem 1 reduces to the following result.
\begin{theorem}
Let it be $D_0(s) = (n - 1)sF'_0(s) - aF_0(s)$.\\
1) If $D_0(n^{\frac{1}{a}}) > 0, \text{ then } Q > 0$. \\
2) If $D_0(1 + \frac{n - 1}{a}) <0, \text{ then }  Q = 0$.
\end{theorem}

\noindent If $a = 1$, then the conditions of points 1 and 2 stick together, and we get the result from \cite{lifsh}. If $a > 1$, then Theorem 2 is a new result.

\section{Proof of the main result}
\subsection{Moment generating functions and their evolution}

Let $F_n(s)$ be the moment generating function of $X_n$. \\
Let us rewrite the recurrence relation \eqref{1} in terms of generating functions:

\begin{eqnarray}
\nonumber
 F_{n+1}(s) &=& \frac{G(F_n(s)) - \sum\limits_{p = 0}^{a - 1}\frac{s^p}{p!}(G(F_n)^{(p)}(0))}{s^a} + \sum_{p = 0}^{a - 1}\frac{G(F_n)^{(p)}(0)}{p!}
\\
\label{3}
&=& \frac{G(F_n(s))}{s^a} + \sum_{p = 0}^{a - 1}(G(F_n))^{(p)}(0)\cfrac{1}{p!}\left(1 - \frac{1}{s^{a-p}}\right).
\end{eqnarray}

By differentiating equality \eqref{3}  we get
\begin{eqnarray}
\label{4}
F'_{n+1}(s) = \frac{G'(F_n(s))F'_n(s)}{s^a} - a\frac{G(F_n(s))}{s^{a+1}} + \sum_{p = 0}^{a - 1}G(F_n)^{(p)}(0) \cdot \cfrac{a-p}{p!} \cdot \cfrac{1}{s^{a-p+1}}\text{.} 
\end{eqnarray}
\noindent We will use these formulas a lot.

\subsection{Proof of the first point of Theorem 1}
\subsubsection*{}
\label{sec:3rd}

The proof relies on two lemmas. First, we will prove the result using them, and the proof
of the lemmas themselves will be presented below.

\begin{lemma}
If the assumption of the first statement of the theorem is true, then there exists  $1 < s < (\mathbb{E}N)^\frac{1}{a}$, such that

$$(\mathbb{E}N - 1)\mathbb{E}(X_ns^{X_n}) - a\mathbb{E}(s^{X_n}) \rightarrow \infty \text{ as } n \rightarrow \infty.$$

\end{lemma}
\noindent In terms of generating functions, Lemma 1  can be represented as follows: there is $1 < s < (\mathbb{E}N)^\frac{1}{a}$ such that

$$(\mathbb{E}N - 1)sF'_n(s) - aF_n(s) \rightarrow \infty \text{ as } n \rightarrow \infty .$$

\begin{lemma}
    If $Q = 0$, then $\sup\limits_{n \geq 0}\mathbb{E}(X_ns^{X_n}) < \infty$ with $0 < s < (\mathbb{E}N)^\frac{1}{a}$.
\end{lemma}

\noindent Now let us apply these lemmas to prove the statement.

\noindent \\Let $a\mathbb{E}(((\mathbb{E}N)^\frac{1}{a})^{X_0}) < (\mathbb{E}N - 1)\mathbb{E}(X_0((\mathbb{E}N)^\frac{1}{a})^{X_0})$; we will show that $Q > 0 $.

\noindent \\By Lemma 1, there is $1 < s < (\mathbb{E}N)^\frac{1}{a}$, such that $ (\mathbb{E}N - 1)\mathbb{E}(X_ns^{X_n} ) - a\mathbb{E}(s^{X_n}) \rightarrow \infty, \text{ as $n$} \rightarrow \infty ,$  therefore $(\mathbb{E}N - 1)\mathbb{E}(X_ns^{X_n}) \rightarrow \infty \text{ as $n$} \rightarrow \infty$.

\noindent \\But if $Q = 0$, then by Lemma 2 $\sup\limits_{n \geq 0}\mathbb{E}(X_n s^{X_n}) < \infty$. Hence $Q > 0$.

\noindent Now we will provide proof of the lemmas.

\subsubsection*{Proof of Lemma 1}

By \eqref{4}

$$s(\mathbb{E}N - 1)F'_{n+1}(s) = \frac{(\mathbb{E}N - 1)sG'(F_n(s))F'_n(s)}{s^a} - a(\mathbb{E}N - 1)\frac{G(F_n(s))}{s^{a}} + \sum_{p = 0}^{a - 1}(\mathbb{E}N - 1)G(F_n)^{(p)}(0) \cdot \cfrac{a-p}{p!} \cdot \cfrac{1}{s^{a-p}}\text{.}$$

\noindent Let us give the lower bound for $G'(F_n(s))$.

\noindent Let us prove the inequality: for any $v \geq 1$, $vG'(v) \geq \mathbb{E}N\cdot G(v)$. In terms of generating functions, it means: $\mathbb{E}(Nv^{N}) \geq \mathbb{E}N \cdot E(v^N)$.

\noindent For any independent copies $N_1, N_2$ of the random variable $N$ the following is true:

\noindent If $(N_1 - N_2)(v^{N_1} - v^{N_2}) \geq 0$, then $0 \leq \mathbb{E}[(N_1 - N_2)(v^{N_1} - v^{N_2})] = 2(\mathbb{E}(Nv^{N}) - \mathbb{E}N\mathbb{E}v^{N})$.

\noindent We note that for any $v \geq 1$:   $(N_1 - N_2)(v^{N_1} - v^{N_2}) \geq 0$.
\newline
\newline
\noindent Thus, for any $v \geq 1$, $vG'(v) \geq \mathbb{E}N \cdot G(v)$. Therefore, $F_n(s)G'(F_n(s)) \geq \mathbb{E}N \cdot G(F_n(s))$.
 Let us substitute this relation into the previous inequality:

 $$s(\mathbb{E}N - 1)F'_{n+1}(s) \geq \frac{\mathbb{E}N(\mathbb{E}N - 1)sG(F_n(s))F'_n(s)}{F_n(s)s^a} - a(\mathbb{E}N - 1)\frac{G(F_n(s))}{s^{a}} + \sum_{p = 0}^{a - 1}(\mathbb{E}N - 1)G(F_n)^{(p)}(0) \cdot \cfrac{a-p}{p!} \cdot \cfrac{1}{s^{a-p}}\text{.}$$
Therefore,
\label{5}
$$s(\mathbb{E}N - 1)F'_{n+1}(s) - aF_{n + 1}(s)$$ $$ \geq \frac{\mathbb{E}N \cdot G(F_n(s))}{F_n(s)s^a}\cdot\left(s(\mathbb{E}(N) - 1)F'_n(s) - aF_n(s)\right)  + \sum_{p = 0}^{a - 1}G(F_n)^{(p)}(0) \cdot \cfrac{1}{p!} \cdot \cfrac{(a-p)(\mathbb{E}N - 1) + a - as^{a-p}}{s^{a-p}}.\eqno{(5)}$$

\noindent Let us prove that $(a - p)(\mathbb{E}N - 1) + a - a s^{a - p} \geq0 $ for $p = 0, 1,\dots, a - 1$ and $s \leq \mathbb{E}N^{\frac{1}{a}}$.

\noindent \\Rewrite this as:
\label{6}
$$y(\mathbb{E}N - 1)+a \geq as^y \text{,   where }y=a-p > 0.\eqno{(6)}$$

\noindent Let us introduce two functions $f_1(y)= a(\mathbb{E}N)^{\frac{y}{a}} \text{ and } 
f_2(y)= y(\mathbb{E}(N) - 1)+a$.

\noindent These functions are equal at $y=0$ and $y=a$ and the function $f_1$ is concave, while $f_2$ is linear.

\noindent Therefore,
$$f_2(y) \geq f_1(y), \quad 0 \leq y \leq a,$$
namely,
 $$y(\mathbb{E}N - 1)+a \geq a (\mathbb{E}N)^{\frac{y}{a}} \geq as^y.$$

\noindent We have proved inequality \hyperref[6]{(6)}, hence the sum in \hyperref[5]{(5)} is non-negative.

\noindent We conclude that
 
$$(\mathbb{E}N - 1)sF'_{n + 1}(s) - aF_{n + 1}(s) \geq \frac{\mathbb{E}N \cdot G(F_n(s))}{F_n(s)s^a}[(\mathbb{E}N - 1)sF'_n (s) - aF_n(s)].$$

\noindent Let $G(v) = \sum\limits_{k = 1}^{\infty}a_{k}v^{k}$, where $a_{k} \geq 0$ ($N > 0$). So $\frac{G(v)}{v} = \sum\limits_{k = 1}^{\infty}a_{k}v^{k - 1} \geq \sum\limits_{k = 1}^{\infty}a_{k}1^{k-1} = 1$.

\noindent This means that
 
$$(\mathbb{E}N - 1)sF'_{n + 1}(s) - aF_{n + 1}(s) \geq \left(\frac{\mathbb{E}N}{s^a}\right)^n [(\mathbb{E}N - 1)sF'_0(s) - aF_0(s)].$$

\noindent Let $p = \frac{\mathbb{E}N}{s^a} > 1$, $C_1 = (\mathbb{E}N - 1)sF'_0(s) - aF_0(s)$, then 
$$(\mathbb{E}N - 1)sF'_{n + 1}(s) - aF_{n + 1}(s) \geq p^nC_1,$$

\noindent which goes to $+\infty$ , if $C_1 > 0 $. As by assumption $D_0(\mathbb{E}(N)^{\frac{1}{a}}, \mathbb{E}N) > 0$, hence, for $s$ close to $\mathbb{E}N^{\frac{1}{a}}$ we have $D_0(s, \mathbb{E}N) > 0$.

\subsubsection*{Proof of Lemma 2}

\noindent \\Let us fix $k \geq 1$ and $n \geq 0.$ Consider the following inequality
$$X_{n + k} \geq \sum_{i = 1}^{T_{k}}\mathds{1}_{X_n^{(i)} \geq ak + 1}, $$

\noindent where $X_n^{(i)}, i \geq 1,$ are independent copies of $X_n$, and $T_{k}$ is a random variable equal to the number of ancestors, when building our counting tree, at depth $k$. It follows, from independence and properties of mathematical expectation that $\mathbb{E}T_{k} = (\mathbb{E}N)^{k}$.

\noindent \\Therefore,
$$\mathbb{E}(X_{n + k}) \geq (\mathbb{E}N)^{k}\mathbb{P}(X_n \geq ak + 1).$$

\noindent On the other hand,  because $Q = 0$, we may use \hyperref[2]{(2)}, hence $\mathbb{E}(X_{n + k}) \leq \frac{a}{\mathbb{E}N - 1}$ for any $n \geq 0, k \geq 1$.

\noindent Therefore,
$$\mathbb{P}(X_n \geq ak + 1) \leq \frac{a}{(\mathbb{E}N - 1)(\mathbb{E}N)^k}.$$

\noindent Let us sum it up and get the desired inequality.
$$\mathbb{E}(X_ns^{X_n}) = \sum_{k = 0}^{\infty}\mathbb{P}(X_n = k)ks^k \leq \sum_{k = 1}^{\infty}\mathbb{P}\Big(ak + 1 \geq X_n > a(k - 1) + 1\Big)(ak + 1)s^{ak + 1} $$  $$ \leq \sum_{k = 1}^{\infty}\mathbb{P}(X_n \geq a(k - 1) + 1)(ak + 1)s^{ak + 1} \leq \sum_{k = 1}^{\infty}\frac{a}{(\mathbb{E}N - 1)(\mathbb{E}N)^{k - 1}}(ak + 1)s^{ak + 1}$$ $$ = s\left(\sum_{k = 1}^{\infty}\frac{a}{(\mathbb{E}N - 1)(\mathbb{E}N)^{k - 1}}s^{ak + 1}\right)' = C_{s, a, \mathbb{E}N} < \infty,$$

\noindent where $C_{s, a, \mathbb{E}N}$ is some constant that depends on $s, \mathbb{E}N$ and $a$.

\noindent Thus,
$$\mathbb{E}(X_ns^{X_n})  \leq C_{s, a, \mathbb{E}N}.$$

\subsection{Proof of the second point of Theorem 1}
\subsubsection*{}
\label{sec:4th}
Let us define the sequence of functions $D_n$:
$$D_n(s) = (M - 1)sF'_n(s) - aF_n(s).$$
As in Section 3, we formulate two lemmas, which will be proved below.
 \begin{lemma}
     $$\text{For any } s \geq 1+ \cfrac{M - 1}{a}, \ \text{and } n\geq 0 \text{ it is true that }  D_{n+1}(s) \leq  \frac{MG(F_n(s))}{F_n(s)s^a}D_n(s).$$  
 \end{lemma}
 \begin{lemma}
    If $s > 1$, then $$\mathbb{E}(X_ns^{X_n}) \geq \mathbb{E}(X_n)\mathbb{E}(s^{X_n}).$$
 \end{lemma}

\noindent Lemma 3 implies that if $D_0(s_0) < 0$ for some $s_0 \geq 1 + \cfrac{M - 1}{a}$, 
then for any $n > 0$ inequality $D_n(s_0) < 0$ holds.

\noindent Let us express $D_n$ in terms of mathematical expectation:
$$D_n(s) = (M - 1)sF'_n(s) - aF_n(s) = (M - 1)\mathbb{E}(X_ns^{X_n}) - a\mathbb{E}(s^{X_n}) \leq 0,$$
Therefore,
$$a\mathbb{E}(s^{X_n}) \geq (M - 1)\mathbb{E}(X_ns^{X_n}).$$

\noindent By Lemma 4 we have:
$$a\mathbb{E}(s^{X_n}) \geq  (M - 1)\mathbb{E}(X_n)\mathbb{E}(s^{X_n}).$$

\noindent We get that $\mathbb{E}(X_n) \leq \frac{a}{M - 1}$, which means that $Q = 0$.

\subsubsection*{Proof of Lemma 3}

\ \newline
\noindent For $p = 0, 1, \dots, a - 1$ and $s \geq 1 + \frac{M - 1}{a}$ we will prove the following:
$$(a-p)(M - 1) + a - a s^{a-p} \leq 0.$$
\noindent Let us denote $y = a - p \geq 1$. and rewrite the required fact as
$$y(M - 1)+a \leq a s^y.$$

\noindent \hypertarget{ineq:1}{Let us use Bernoulli’s inequality:}

$$(1 + x)^n \geq 1 + nx, \ \text{for any } n \in \mathbb{N}, x > 0.$$

\noindent We note that $as^y = a(1+(s-1))^y.$

\noindent By assumption $y \geq 1, s \geq 1+\frac{1}{a}; \text{ it means } s - 1 \geq \frac{M - 1}{a} > 0 $. Therefore, we can apply Bernoulli's inequality to $ x = s - 1$ and $n = y$. It gives us:
$$(1+(s-1))^y \geq 1+(s-1)y \geq 1 + \frac{y(M - 1)}{a}.$$
Hence,
$$as^y \geq a\left(1 + \frac{y(M - 1)}{a}\right) = (M - 1)y + a.$$

\noindent Let us prove the auxiliary inequality:

$$vG'(v) \leq MG(v), \text{ namely } F_n(s)G'(F_n(s)) \leq MG(F_n(s)).$$

\noindent As $\mathbb{E}N \leq M$ is true,  $vG'(v) \leq \mathbb{E}(Nv^{N}) \leq
\mathbb{E}N \cdot \mathbb{E}v^{N} \leq M \cdot \mathbb{E}v^{N} = MG(v)$ is true.

\noindent It follows that
$$s(\mathbb{E}(N) - 1)F'_{n+1}(s) - aF_{n + 1}(s)$$ $$\leq \frac{MG(F_n(s))}{F_n(s)s^a}\cdot\left(s(M - 1)F'_n(s) - aF_n(s)\right)  + \sum_{p = 0}^{a - 1}G(F_n)^{(p)}(0) \cdot \cfrac{1}{p!} \cdot \cfrac{(a-p)(M - 1) + a - as^{a-p}}{s^{a-p}}.\eqno{(7)}$$

\noindent All terms in the corresponding sum are non-positive, because as proven, $(a-p)(M-1) + a - as^{a-p} \leq 0$.

\noindent Therefore,

$$D_{n+1}(s) \leq  \frac{MG(F_n(s))}{F_n(s)s^a}D_n(s).$$

\subsubsection*{Proof of Lemma 4}

\noindent \\Let $Y_1$, $Y_2$ be independent copies $X_n$. 

\noindent \\At $s \geq 1$ we have:

$$(Y_1 - Y_2)(s^{Y_1} - s^{Y_2}) \geq 0.$$

\noindent \\Therefore,
$$\mathbb{E}[(Y_1 - Y_2)(s^{Y_1} - s^{Y_2})] \geq 0.$$

\noindent $Y_1$ and $Y_2$ are independent, hence $Y_1$ and $s^{Y_2}$ are independent.

\noindent \\Therefore,

$$\mathbb{E}(Y_1s^{Y_1}) + \mathbb{E}(Y_2s^{Y_2}) - \mathbb{E}Y_1\mathbb{E}s^{Y_2} - \mathbb{E}Y_2\mathbb{E}s^{Y_1} \geq 0.$$

\noindent $Y_1$, $Y_2$ are independent copies $X_n$, hence

$$2\mathbb{E}(X_ns^{X_n}) - 2\mathbb{E}X_n\mathbb{E}s^{X_n} \geq 0,$$

$$\mathbb{E}(X_ns^{X_n}) - \mathbb{E}X_n\mathbb{E}s^{X_n} \geq 0.$$

\noindent The authors are grateful to M. Lifshits for setting the problem and useful advice. \\

St. Petersburg State University,

Universitetskaya emb. 7/9, 191023, St. Petersburg, Russia

E-mail : Kotann2710@mail.ru

E-mail : alex.lotnikov@gmail.com
\end{document}